\documentclass{elsart}
\usepackage{amssymb}
\usepackage{amsmath}
\usepackage{amscd}
\usepackage{graphics}
\usepackage{graphicx}
\newtheorem{theo}{Theorem}

\newcommand{\hide}[1]{}

\renewcommand{\emph}{\textbf}

\begin{document}
\begin{frontmatter}
\title{Symbolic lumping of some\\ catenary, mamillary and circular\\
compartmental systems}
\author[egri]{Edith Egri\thanksref{ize}}
\thanks[ize]{The research upon which the present paper is based has been started during a visit of EE
partially supported by the National Scientific Foundation, Hungary,
under No. T 047132.} \ead{egriedit@yahoo.com}
\address[egri]{Babe\c{s}--Bolyai University, Department of Differential Equations,
Cluj Napoca, Str. M. Kog\u{a}lniceanu, nr.1, 3400, ROMANIA}

\author[toth]{J\'anos T\'oth\thanksref{ize}}
\ead{jtoth@math.bme.hu}
\address[toth]{Department of Mathematical Analysis, Budapest University of Technology and Economics,
Budapest, H-1111 Egry J. u. 1., HUNGARY}

\author[brochot]{C\'eline Brochot and Frederic Yves Bois}
\ead{cbrochot@ineris.fr} \ead{frederic.bois@ineris.fr}
\address[brochot]{INERIS,
Institut National de l'Environnement Industriel et des Risques,
Unit\'e de Toxicologie Exp\'erimentale, Parc Alata BP2, 60550
Verneuil En Halatte, FRANCE}

\begin{abstract}
Some of the most important compartmental systems,
such as irreversible catenary, mamillary and circular systems
are symbolically simplified by the method of exact linear lumping.
A few symbolically unmanageable systems are numerically lumped.
Transformation of the qualitative properties under lumping are also traced.
\end{abstract}

\end{frontmatter}
\tableofcontents
\section{Introduction}
Compartmental systems are mathematical systems that are frequently
used in biology and mathematics. Also a subclass of the class of
chemical processes can be modeled as compartmental systems. A
compartmental system consists of several compartments with more or
less homogeneous amounts of material. The compartments interact by
processes of transport and diffusion. The dynamics of a
compartmental system is derived from mass balance considerations.

The mathematical theory of compartmental systems is of major
importance: it is the bread-and-butter of analysis for medical
researchers, pharmacokineticists, physiologists, ecologists,
economists as well as other researchers\\ \cite{cob76},
\cite{dist87}, \cite{cob87}, \cite{jack99}, \cite{nestorov98}.


Sometimes it is useful to reduce a model to get a new one with a
lower dimension. The technique's name is lumping, i.e. reduction
of the number of variables by grouping them via a linear or
nonlinear function.


The objective of model reduction methods is to obtain a model that
can describe the response of the original model accurately and
efficiently (cf. \cite{wilkinson08}).

Our aim here is to give explicitly possible lumped compartmental
systems in a few important classes, mainly of symmetric structure such as: mamillary models,
catenary models and circular models. Some classes can be treated in
full generality, some only under restrictions on the parameters.

We also show how to lump systems which are only numerically
lumpable.

The structure of our paper is as follows. In Section 2 the formal
definitions of reactions, compartmental systems, induced kinetic
differential equations and that of exact linear lumping are given.
Next, our symbolic results are presented. Section 4 shows a few
examples which had to be treated numerically. Finally, the results
are discussed and further goals are set. We mention that the present
work is a continuation of a few simple statements in \cite{btb} on
the symbolic lumping of a general two compartment model.
\section{Fundamental definitions}
\subsection{Reaction mechanism, compartmental system}

A \textit{chemical reaction mechanism} is a set of elementary
reactions. Formally, it is a system
$<\mathcal{M,R},\alpha,\beta>$, where
\begin{enumerate}
\item $\mathcal{M}$ and $\mathcal{R}$ are sets with $M$ and $R$
elements ($M,R\in \mathbb{N}$), $\mathcal{R}=\{1,2,\ldots,R\}$ and
$\mathcal{M}=\{\mathcal{X}_1,\mathcal{X}_2,\ldots,\mathcal{X}_M\},$
\item $\alpha$ and $\beta$ are matrices with non-negative
integers, whose names are \textit{stoichiometric
coefficients}\index{stoichiometry coefficients}, and for which
\begin{enumerate}
\item for all $r\in \mathcal{R}$, $\alpha(.,r)\neq\beta(.,r),$
\item if $\alpha(.,r)=\alpha(.,r^{'})$ and
$\beta(.,r)=\beta(.,r^{'}),$ then $r=r^{'},$ \item for all $m\in
\mathcal{M}$ there exists $r\in \mathcal{R}$ such that either
$\alpha(m,r)\neq 0$ or $\beta(m,r)\neq 0$ holds.
\end{enumerate}
\end{enumerate}

This mechanism can be represented in the form
\begin{equation} \label{formal_mechanism}
\sum_{m=1}^{M}\alpha(m,r)\mathcal{X}_{m}\longrightarrow
\sum_{m=1}^{M}\beta(m,r)\mathcal{X}_{m} \qquad (r\in \mathcal{R}).
\end{equation}
The entities on the two sides of the arrow are the \textit{reactant}
and \textit{product complexes}, respectively.

The number $\max\{\sum_{m=1}^{M}\alpha(m,r),r\in\mathcal{R}\}$ is
said to be the \textit{order} of the reaction; thus, \textit{first
order reaction}s are obtained if \(\forall r\in\mathcal{R}\)
\(\sum_{m=1}^{M}\alpha(m,r)\le 1.\) If in a first order reaction it
is also true that the length \(\sum_{m=1}^{M}\beta(m,r)\) of the
product complexes is also less than or equal to 1, then one has a
\textit{compartmental system}. These formal mechanisms are of great
practical importance, and are applied in many areas as mentioned in
the introduction.

Thus, a compartmental system is a reaction mechanism in which the
length of all the complexes is not more than one. In this case we
only have reaction steps of the type
$\mathcal{X}_m\rightarrow\mathcal{X}_p,\,\,
\mathcal{X}_m\rightarrow\mathcal{O},\,\,
\mathcal{O}\rightarrow\mathcal{X}_m \,\,(m,p\in\mathcal{M}),$
where $\mathcal{O}$ is the empty complex.

A \textit{generalized compartmental system} is a reaction in which
all the complexes contain a single species, and all the species
are contained in a single complex, i.e. it is a reaction
consisting of elementary reactions of three types
\begin{equation}\label{gencomp}
y_m\mathcal{X}_m\rightarrow y_p\mathcal{X}_p,\quad
y_m\mathcal{X}_m\rightarrow\mathcal{O},\quad \mathcal{O}\rightarrow
y_m\mathcal{X}_m, \quad (m,p\in\mathcal{M}),
\end{equation}
and \(\mathcal{X}_m\) is the constituent of a single complex only.

A generalized compartmental system with no inflow and with some
outflow is \textit{strictly half-open}, while it is
\textit{strictly open} if it contains inflows and possibly
outflows.

Reaction \eqref{formal_mechanism} is said to be
\textit{mass-conserving} if there exist positive numbers $\rho(1),
\rho(2)$, $\ldots, \rho(N)$ such that for all elementary reactions
\begin{equation}
\sum_{m=1}^{M}\alpha(m,r)\rho(m)= \sum_{m=1}^{M}\beta(m,r)\rho(m)
\end{equation}
holds. If the atomic structure of the species are not known, it is
not trivial to decide whether a reaction is mass-conserving or not
\cite{deak}, \cite{schuster}.

A generalized compartmental system is mass-conserving if and only if
it is closed: the empty complex is not present.

\subsection{Induced kinetic differential equations}
The usual continuous time, continuous state deterministic model
(or, induced kinetic differential equation) of reaction
\eqref{formal_mechanism} describing the time evolution of the
concentrations $c_m$ is the polynomial differential equation
\begin{equation}
\dot{c}_m=\sum_{r=1}^{R}(\beta(m,r)-\alpha(m,r))k_r\prod_{p=1}^M
c_p^{\alpha(p,r)},
\end{equation}
where $k_r$ denotes the rate coefficient, for all $r\in\mathcal{R}.$

The induced kinetic differential equation of a first order reaction
is of the form
\begin{equation}\label{linear}
\dot{c}=Ac+b
\end{equation}
with
\begin{equation}\label{req1}
a_{mp}\geq 0 \quad (m\neq p) \quad \mbox{and} \quad b_m\geq 0
\quad (m,p\in\mathcal{M}).
\end{equation}
The induced kinetic differential equation of a compartmental system
has an additional property
\begin{equation}\label{req2}
-a_{mm}\geq \sum_{\substack{p=1\\p\neq m}}^Ma_{pm}\qquad (m\in
\mathcal{M}).
\end{equation}
Thus, e.g. there is no compartmental system with the induced kinetic
differential equation $\dot{x}=x$ or with $\dot{x}=-0.5x+y, \quad
\dot{y}=-y+x.$

An easy construction proves that the converse of the above
statement is also true: a linear differential equation
\eqref{linear} fulfilling the requirements \eqref{req1} and
\eqref{req2} can be considered as the induced kinetic differential
equation of a compartmental system.

This statement can be generalized to get our next theorem showing
that if the right hand side of a kinetic differential equation is
the sum of univariate monomials and if all the variables have the
same exponent in all the rows, then -- if an additional condition
is also met and only then -- there exists an inducing generalized
compartmental system to the system of differential equations.

\begin{theo}
There exists an inducing generalized compartmental system of $M$
compartments to the system of differential equations
\begin{equation}\label{comp}
\dot{c}_m=\sum_{p=1}^{M}a_{mp}(c_p)^{y^p}+b_m
\end{equation}
(where for all $m,p\in\mathcal{M}, y^m, y^p\in\mathbb{N}, y^m\neq
y^p, if m\neq p, a_{mp},b_m\in\mathbb{R}$) which is
\begin{enumerate}
\item closed, if and only if $b_m=0, -a_{mm},a_{mp},d_m\in
\mathbb{R}_0^+; a_{mm}=d_my^m,$ \item strictly half-open, if and
only if $b_m=0, -a_{mm},a_{mp},d_m\in \mathbb{R}_0^+; a_{mm}\leq
d_my^m, \exists m, a_{mm}<d_my^m,$ \item strictly open, if and
only if $b_m, -a_{mm},a_{mp},d_m\in \mathbb{R}_0^+; a_{mm}\leq
d_my^m, \exists m \, b_m\in \mathbb{R}^+,$
\end{enumerate}
where throughout $$m,p\in\mathcal{M}, m\neq p, d_m:=
-\sum_{p=1}^M a_{pm}/y^m.$$
\end{theo}
\textbf{Proof.}
\begin{enumerate}
\item [\textbf{A)}]The induced kinetic differential equation of
\eqref{gencomp} is
 \begin{equation*}
\begin{split}
\dot{c}_m=&-y^m(c_m)^{y^m}\sum_p k_{pm}+y^m\sum_p k_{mp}(c_p)^{y^p},\\
\dot{c}_m=&-y^m(c_m)^{y^m}\sum_p k_{pm}+y^m\sum_p k_{mp}(c_p)^{y^p},\quad (\exists\, k_{0m}\in\mathbb{R}^+)\\
\dot{c}_m=&-y^m(c_m)^{y^m}\sum_p k_{pm}+y^m\sum_p k_{mp}(c_p)^{y^p}+k_{m0},\quad (\exists\, k_{m0}\in\mathbb{R}^+)\\
& m\in\{1,2,\ldots,\mathcal{M}\};\, k_{mp}\in\mathbb{R}_0^+;\,y^p\in\mathbb{N}; p\in\{0,1\ldots,\mathcal{M}\}
\end{split}
\end{equation*}
Comparing the coefficients we get the only if part of the Theorem.
 \item[\textbf{B)}] Given \eqref{comp} we construct a generalized
 compartmental system \eqref{comp} as its induced kinetic
 differential equation:
\begin{equation}\label{const}
y^p\mathcal{X}_p\stackrel{a_{mp}/y^m}{\longrightarrow} y^m\mathcal{X}_m,\quad
y^p\mathcal{X}_p\stackrel{d_p}{\longrightarrow}\mathcal{O},\quad
\mathcal{O}\stackrel{b_m/y^m}{\longrightarrow} y^m\mathcal{X}_m,
\end{equation}
$(m,p\in\{1,2,\ldots,\mathcal{M}\},m\neq p).$
\end{enumerate}
Reaction \eqref{const} induces closed, strictly half-open or
strictly open reactions, respectively.

\subsection{Exact linear lumping}

A special class of lumping is exact linear lumping.

A system $\dot{c}=f\circ c,$ with $f,c$ $n$-vectors can be
\textit{exactly lumped} by an $\hat{n}\times n$ real constant
matrix $Q$ ($\hat{n}<n$), called \textit{lumping matrix}, if for
$\hat{c}=Qc$ we can find an $\hat{n}$-function vector $\hat{f}$
such that $ \dot{\hat{c}}=\hat{f}\circ \hat{c}.$

Not every system is exactly lumpable. A sufficient and necessary
condition for the existence of exact lumping is $Qf(c)=Qf
(\overline{Q}Qc)$, where $\overline{Q}$ denotes any of the
generalized inverses of $Q$, i.e. $Q\overline{Q}=I_{\hat{n}},$ and
$I_{\hat{n}}$ is the $\hat{n}\times \hat{n}$ identity matrix
\cite{lirab89}.

This condition is equivalent to the requirement that the rows of
matrix $Q$ span an invariant subspace of ${f'}^{\top}(c)$ for all
$c$, where ${f'}^{\top}(c)$ denotes the transpose of the Jacobian of
$f$ at $c$. Therefore, in order to determine lumping matrices $Q$ we
need to determine the fixed ${f'}^{\top}(c)$-invariant subspaces\\
\cite{gohberg}.

In the case of linear differential equation \eqref{linear}, the
Jacobian matrix is just $A$, and then ${f'}^{\top}(c)=A^{\top}.$ In
this situation, fixed invariant subspaces exist, they are spanned by
eigenvectors, and they correspond to (constant) eigenvalues. So, a
linear system is always exactly lumpable and any
${f'}^{\top}(c)$-invariant subspaces will give a lumping matrix. In
this case therefore, we have to calculate the eigenvectors of
$A^{\top}.$

We mention here that, if $Q$ is an $\hat{n}\times n$ lumping matrix
and $P$ a nonsingular matrix of dimension $\hat{n},$ then $PQ$ is
also a lumping matrix.

It is not true that a given system can be lumped arbitrarily. For
example
\begin{equation}\label{i1i2}
\begin{split}
&S\stackrel{k_1}{\rightarrow} I_1\stackrel{k_2}{\rightarrow} P,\\
&S\stackrel{k_3}{\rightarrow} I_2\stackrel{k_4}{\rightarrow} P,\\
\end{split}
\end{equation}
cannot lead to a lumped system of the type
$S\stackrel{K_1}{\rightarrow} I\stackrel{K_2}{\rightarrow}P,$ for
$I:=I_1+I_2,$ except in the very special case $k_2=k_4$, contrary to\\
\cite{Conzelmann}.

The next question is whether the lumped system can have an
interpretation in terms of reactions (or more specially, in terms of
compartments), i.e. is the lumped system kinetic? To formulate this
criterion we use the notion of the generalized inverse matrix
\cite{rao73}.

Farkas \cite{fgy99} gave a sufficient and necessary condition
under which certain lumping schemes preserve the kinetic structure
of the original system: \textit{A nonnegative lumping matrix leads
to a kinetic differential equation if and only if it has a
nonnegative generalized inverse.}

For the absence of a nonnegative generalized inverse he proved
the following result: \textit{A nonnegative matrix has no
nonnegative generalized inverse if and only if it has a row such
that in the column of each positive entry there exists another
positive entry.}

\section{Symbolic results}
\subsection{Chains}
In a chain or catenary system the $M$ compartments are arranged in a
linear array such that every compartment exchanges  material only
with its immediate neighbors and the possible steps are indicated by
nonnegative reaction rates. The coefficient matrix for a catenary
system has nonzero entries only in the main diagonal and the first
sub-diagonal and in the first super-diagonal. The latter case holds
only if it is reversible (bidirectional).

\subsubsection{Irreversible chains}
Let us consider a compartmental system, such as the one in
Fig.~\ref{IrrCat}, i.e. a chain with unidirectional steps.

\begin{figure}[h!]
\begin{center}
\includegraphics[totalheight=0.4in]{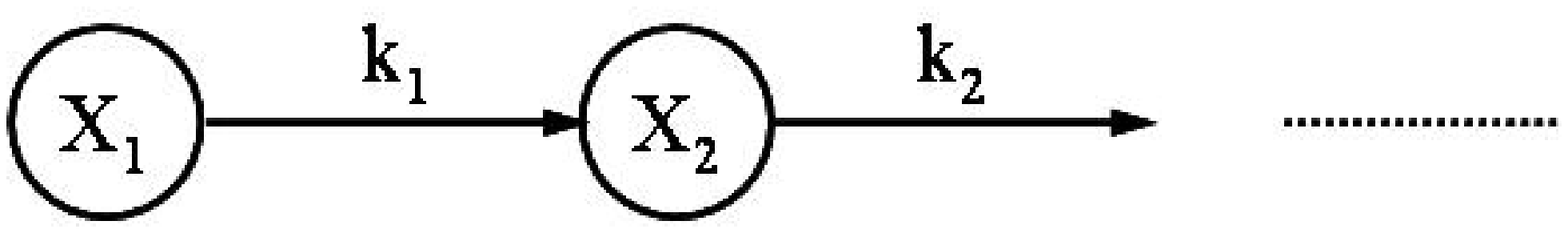}
\end{center}
\caption{Irreversible catenary system}\label{IrrCat}
\end{figure}

In this case the coefficient matrix $A$ on the right hand side of
\eqref{linear} takes the form
\[
 \left[
  \begin{array}{ccccccc}
  -k_1    & 0      & 0      & \ldots & 0          & 0        & 0\\
   k_1    & -k_2   & 0      & \ldots & 0          & 0        & 0\\
   0      &  k_2   & -k_3   & \ldots & 0          & 0        & 0\\
   \vdots & \vdots & \vdots &        & \vdots     & \vdots   & \vdots\\
   0      & 0      & 0      & \ldots & k_{M-2}    & -k_{M-1} & 0\\
   0      & 0      & 0      & \ldots & 0          & k_{M-1}   & 0
  \end{array}
 \right].
\]

The eigenvalues of the transpose of this triangular matrix are
obviously the elements on the diagonal: $-k_1, -k_2, -k_3, \ldots,
-k_{M-1}$ and $0$ (and are the same as the eigenvalues of the original matrix). So, the corresponding eigenvectors can be found
easily, and they take the form:
\[
  \begin{array}{cccccc}
  \left[ 1      \right.& 0                    & 0                                  & \ldots & 0 &\left.  0\right]     \\
  \left[ 1      \right.& \frac{k_1-k_2}{k_1}  & 0                                  & \ldots & 0 &\left.  0\right]     \\
  \left[ 1      \right.& \frac{k_1-k_3}{k_1}  & \frac{(k_1-k_3)(k_2-k_3)}{k_1k_2}  & \ldots & 0  &\left.  0\right]    \\
  \vdots               & \vdots               & \vdots                             & \vdots & \vdots & \vdots \\
  \left[1       \right.& \frac{k_1-k_{M-1}}{k_1}&\frac{(k_1-k_{M-1})(k_2-k_{M-1})}{k_1k_2}&\ldots &\frac{(k_1-k_{M-1})(k_2-k_{M-1})\cdots(k_{M-2}-k_{M-1})}{k_1k_2\ldots k_{M-1}}   &\left.  0\right]     \\
  \left[1       \right.& 1                    & 1                                  & \ldots & 1 &\left.  1\right]     \\
  \end{array}
 \]

\hide{
\begin{displaymath}
[1,0,0,\ldots, 0],[1, -\frac{k_2-k_1}{k_1}, 0, \ldots, 0], [1,
-\frac{k_3-k_1}{k_1}, \frac{(k_3-k_2)(k_3-k_1)}{k_1k_2}, \ldots,
0],\ldots
\end{displaymath}
\begin{displaymath}
[1,-\frac{k_M-k_1}{k_1},\frac{(k_M-k_2)(k_M-k_1)}{k_1k_2},
-\frac{(k_M-k_3)(k_M-k_2)(k_M-k_1)}{k_1k_2k_3}, \ldots, 0],
[1,1,1,\ldots,1].
\end{displaymath}
}

(Here we only consider the robust case when all the reaction rate
coefficients are different. Then, the above eigenvectors are
independent.)

If we do not neglect inflows and outflows in  a catenary system,
the principal diagonal of matrix $A$ will change, i.e. instead of
$-k_i$ we will have $-k_i-\mu_i$ in the first $M-1$ places (where
$\mu_i$ denotes the outflow coefficient for the species $X_i$),
and $-\mu_{M}$ in the last one, instead of $0$. The transpose of
the modified matrix has the following eigenvectors:
\[
  \begin{array}{lcccc}
  \left[ 1      \right.& 0                    & 0                                  & \ldots &\left. 0\right]     \\
  \left[ \frac{k_1}{k_1-k_2+\mu_1-\mu_2}      \right.& 1  & 0                                  & \ldots &\left.  0\right]     \\
  \left[ \frac{k_1k_2}{(k_1-k_3+\mu_1-\mu_3)(k_2-k_3+\mu_2-\mu_3)}      \right.& \frac{k_2}{k_2-k_3+\mu_2-\mu_3}  & 1  & \ldots &\left.  0\right]     \\
  \vdots  & \vdots               & \vdots                             & \vdots & \vdots \\
  \left[1      \right.& 1                    & 1                                  & \ldots &\left.  1\right],     \\
  \end{array}
 \] corresponding to the eigenvalues $-k_1-\mu_1, -k_2-\mu_2, \ldots -k_{M-1}-\mu_{M-1}, -\mu_{M}.$

The graphical representation in this case, when outflows and inflows are incorporated into an irreversible chain, is:

\begin{figure}[h!]
\begin{center}
\includegraphics[totalheight=1.0in]{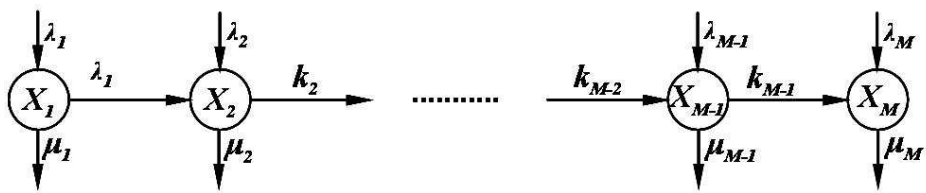}
\end{center}
\caption{Irreversible catenary system with inflows and outflows}\label{IrrCat}
\end{figure}

To get a lumped system for this model, we can take some of the
eigenvectors above to generate several lumping matrices.

For example, let us consider  an irreversible chain with five
compartments. Then the induced kinetic differential equation has the
following coefficient matrix:

\[
 \left[
  \begin{array}{ccccc}
  -k_1-\mu_1   & 0                & 0              &0            &0\\
  k_1          & -k_2-\mu_2       & 0              &0            & 0\\
   0           & k_2              & -k_3-\mu_3     &0            & 0\\
   0           & 0                & k_3            &-k_4-\mu_4   & 0\\
   0           & 0                & 0              &k_4          & -\mu_5
  \end{array}
 \right].
\]

Let us compose $Q$ e.g. putting the  eigenvectors
$\left[\frac{k_1}{k_1-k_2+\mu_1-\mu_2}\quad 1\quad 0\quad 0\quad
0\right]$ and
$\left[\frac{k_1k_2k_3}{(k_1-k_4+\mu_1-\mu_4)(k_2-k_4+\mu_2-\mu_4)(k_3-k_4+\mu_3-\mu_4)}\quad
\frac{k_2k_3}{k_3-k_4+\mu_3-\mu_4}\quad
\frac{k_3}{(k_2-k_4+\mu_2-\mu_4)(k_3-k_4+\mu_3-\mu_4)} \quad 1\quad
0\right]$

into it as rows. Then,
\[Q^T=
 \left[
  \begin{array}{cc}
   \frac{k_1}{k_1-k_2+\mu_1-\mu_2} & \frac{k_1k_2k_3}{(k_1-k_4+\mu_1-\mu_4)(k_2-k_4+\mu_2-\mu_4)(k_3-k_4+\mu_3-\mu_4)}\\
   1    &   \frac{k_2k_3}{(k_2-k_4+\mu_2-\mu_4)(k_3-k_4+\mu_3-\mu_4)}\\
   0    &   \frac{k_3}{k_3-k_4+\mu_3-\mu_4}\\
   0    &   1\\
   0    &   0\\
  \end{array}
 \right].
\]

\hide{
\[Q=
 \left[
  \begin{array}{ccccc}
   \frac{k_1}{k_1-k_2+\mu_1-\mu_2}                           & 1                                & 0   & 0 & 0  \\
   \frac{k_1k_2k_3}{(k_1-k_4+\mu_1-\mu_4)(k_2-k_4+\mu_2-\mu_4)(k_3-k_4+\mu_3-\mu_4)} & \frac{k_1k_2k_3}{(k_1-k_4+\mu_1-\mu_4)(k_2-k_4+\mu_2-\mu_4)(k_3-k_4+\mu_3-\mu_4)}  & \frac{k_3}{k_3-k_4+\mu_3-\mu_4}   & 1 & 0  \\
  \end{array}
 \right].
\]
}

After some calculations, we get the lumped system
$\hat{A}=QA\overline{Q}$, which induces the differential equation
below:
\[\left[
  \begin{array}{c}
  \dot{\hat{x}}_1\\
  \dot{\hat{x}}_2\\
  \end{array}
 \right]=\left[
  \begin{array}{cc}
  -k_2-\mu_2   & 0   \\
  0      & -k_4-\mu_4\\
  \end{array}
 \right]\left[
  \begin{array}{c}
  \hat{x}_1\\
  \hat{x}_2
  \end{array}
 \right],
 \] so we got a new compartmental system with two
 compartments, where
 \begin{equation*}
\begin{split}
\hat{x}_1&=\frac{k_1}{k_1-k_2+\mu_1-\mu_2}x_1+x_2\\
\hat{x}_2&=\frac{k_1k_2k_3}{(k_1-k_4+\mu_1-\mu_4)(k_2-k_4+\mu_2-\mu_4)(k_3-k_4+\mu_3-\mu_4)}x_1+\\
&+\frac{k_2k_3}{(k_2-k_4+\mu_2-\mu_4)(k_3-k_4+\mu_3-\mu_4)}x_2+\frac{k_3}{k_3-k_4+\mu_3-\mu_4}x_3+x_4.
\end{split}
\end{equation*}
The corresponding reaction (actually, a chain with no interaction between the compartments) can be illustrated as follows:
\begin{eqnarray*} \label{example1}
\hat{\mathcal{X}}_1\stackrel{k_2+\mu_2}{\rightarrow}\mathcal{O}
\stackrel{k_4+\mu_4}{\leftarrow}\hat{\mathcal{X}}_2,
\end{eqnarray*}
or it can be the mamillary system
$\hat{\mathcal{X}}_1\rightarrow\hat{\mathcal{X}}_3
\leftarrow\hat{\mathcal{X}}_2,$ with $\hat{\mathcal{X}}_3$ neglected in the induced kinetic differential equation.

\subsubsection{Irreversible chains with nonuniform directions}\label{nonunich}
Let us mention here that the irreversible case with nonuniform
direction of the arrows is simpler than the case of a reversible
chain.

As an example, consider the compartmental system with the following diagram:

\begin{figure}[h!]
\begin{center}
\includegraphics[totalheight=0.45in]{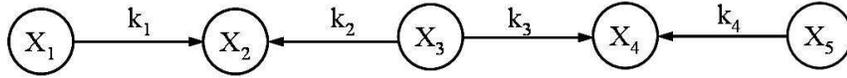}
\end{center}
\caption{An irreversible chain with nonuniform directions}\label{RevCh}
\end{figure}

We can associate to it the kinetic differential equation $\dot{x}=Ax,$ where
\[A=
 \left[
  \begin{array}{ccccc}
  -k_1         & 0                & 0              &0            & 0\\
  k_1          & 0                & k_2            &0            & 0\\
   0           & 0                & -k_2-k_3       &0            & 0\\
   0           & 0                & k_3            &0            & k_4\\
   0           & 0                & 0              &0            & -k_4
  \end{array}
 \right].
\]

The eigenvalues are: $-k_1, -k_2-k_3, -k_4,$ and $0,$ with
multiplicity 2. From the corresponding eigenvectors,
$[1,0,0,0,0],[0,0,1,0,0],[0,0,0,0,1],
\left[\cfrac{k_2+k_3}{k_2},\cfrac{k_2+k_3}{k_2},1,0,0\right]$ and
$\left[-\cfrac{k_3}{k_2},-\cfrac{k_3}{k_2},0,1,1\right]$ we can
determine a lot of lumping matrices. Depending on our choice, the
lumped system can be kinetic or not.

For example, if we take
\[Q=
 \left[
  \begin{array}{ccccc}
   0           & 0                & 0              &0            & 1\\
   0           & 0                & 1              &0            & 0\\
   1           & 0                & 0              &0            & 0
  \end{array}
 \right],
\] then for the lumped system  we get the kinetic differential equation system
\[\left[
  \begin{array}{c}
  \dot{\hat{x}}_1\\
  \dot{\hat{x}}_2\\
  \dot{\hat{x}}_3
  \end{array}
 \right]=\left[
  \begin{array}{ccc}
  -k_4   & 0        & 0\\
  0      & -k_2-k_3 & 0\\
  0      & 0        & -k_1
  \end{array}
 \right]\left[
  \begin{array}{c}
  \hat{x}_1\\
  \hat{x}_2\\
  \hat{x}_3
  \end{array}
 \right],
 \] which can be illustrated via the diagram:
\begin{eqnarray*}
\mathcal{\hat{X}}_1\stackrel{k_4}{\longrightarrow}&&
\mathcal{O}\stackrel{k_1}{\longleftarrow}\mathcal{\hat{X}}_3\\
&&\uparrow \hbox{\scriptsize{$k_2+k_3$}}\\
&&\mathcal{\hat{X}}_2
\end{eqnarray*}

On the other hand, if we take
\[Q=
 \left[
  \begin{array}{ccccc}
   1                & 0                & 0              &0            & 0\\
   0                & 0                & 1              &0            & 0\\
   \cfrac{k_2+k_3}{k_2} & \cfrac{k_2+k_3}{k_2} & 1            &0            & 0
  \end{array}
 \right],
\] this leads
to the matrix
\[\hat{A}=QA\overline{Q}=
 \left[
  \begin{array}{ccc}
   -2k_1                               &-\cfrac{k_1k_2}{k_2+k_3}            & \cfrac{k_1k_2}{k_2+k_3}\\
   -k_2                                &-\cfrac{k_2^2}{k_2+k_3}-k_2-k_3      & \cfrac{k_2^2}{k_2+k_3}\\
   -\cfrac{2k_1(k_2+k_3)}{k_2}-k_2     & -\cfrac{k_1(k_2+k_3)+k_2^2}{k_2+k_3}-k_2-k_3           & \cfrac{k_1(k_2+k_3)+k_2^2}{k_2+k_3}
  \end{array}
 \right].
\]
It can be seen that in this case the positivity conditions relative
to the convenient elements of the matrix are not fulfilled, as
expected in accordance with Lemma 1 in \cite{fgy99}. Consequently,
$\hat{A}$ does not result in a lumped system which has a kinetic
differential equation. The new variables are:
\begin{equation*}
\begin{split}
\hat{x}_1&=x_1\\
\hat{x}_2&=x_3\\
\hat{x}_3&=\frac{k_2+k_3}{k_2}x_1+\frac{k_2+k_3}{k_2}x_2+x_3.
\end{split}
\end{equation*}

\subsubsection{Reversible chains}
\begin{figure}[h!]
\begin{center}
\includegraphics[totalheight=0.6in]{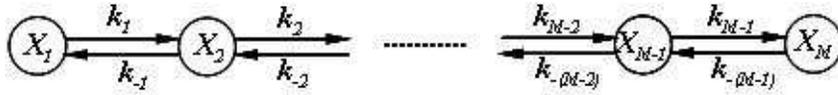}
\end{center}
\caption{Reversible Chain}\label{RevCh}
\end{figure}

To compute the eigenvectors even for a  reversible chain consisting
of only five compartments is unsolvable symbolically. We shall give
a numerical example in section \ref{numrevCh} below.

\subsection{Mamillary systems}

In these systems all the compartments communicate  only with a
central compartment, $X_{M+1}$, and there is no direct communication between the other
compartments. The possible steps are indicated by nonnegative
reaction rates. We shall call $X_{M+1}$ as the \textit{mother
compartment} and all the other compartments will be called
\textit{daughter} or \textit{peripheral} compartments.

\begin{figure}[h!]
\begin{center}
\includegraphics[totalheight=2in]{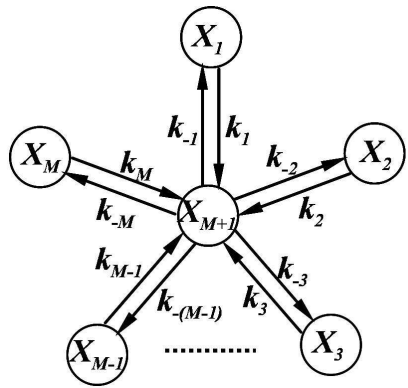}
\end{center}
\caption{Mamillary system}\label{Mam}
\end{figure}

Only the irreversible case can be treated symbolically; a reversible
example will be treated numerically in section \ref{numrevMam}. A
class of reversible mamillary systems with a special structure can
still be treated symbolically, this will be shown in subsection
\ref{classmam}.

\subsubsection{Irreversible mamillary systems}

\textbf{Inward flows}

\begin{figure}[h!]
\begin{center}
\includegraphics[totalheight=2.2in]{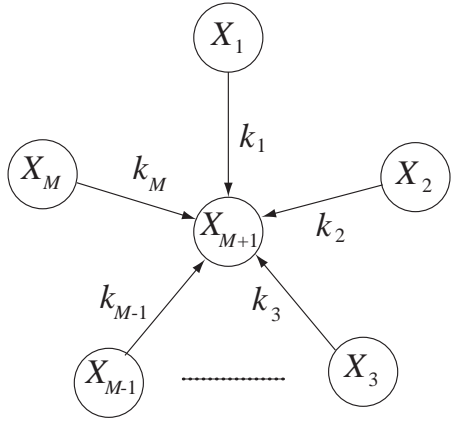}
\end{center}
\caption{A mamillary system with inward flows only}\label{InFlows}
\end{figure}

Let us consider an irreversible mamillary system with inward flows
such as the one in Fig.~\ref{InFlows}.

The coefficient matrix of the reaction rate constants is
\[
 \left[
  \begin{array}{cccccc}
  -k_1   & 0            & \ldots & 0      &0      &0\\
  0      & -k_2         & \ldots & 0      &0      & 0\\
  \vdots & \vdots       & \vdots & \vdots &\vdots &\vdots\\
  0      & 0            & \ldots & 0      &-k_M   & 0\\
  k_1    & k_2          & \ldots & k_{M-1}& k_M   & 0
  \end{array}
 \right].
\]
The eigenvalues of the transpose of this lower triangular matrix are
obviously the elements on the diagonal: $-k_1, -k_2, -k_3, \ldots,
-k_M$ and $0,$ with the corresponding eigenvectors:
\begin{displaymath}
[1,0,\ldots,0,0],[0,1, \ldots,0,0],\ldots, [0,0, \ldots,1,0], [1,1,\ldots,1,1].
\end{displaymath}
Notice that if we denote by $e_i$ the $i$-th element of the standard
basis for $\mathbb{R}^N$, $i\in\{1,2,\ldots, N\}$, then $e_1,
e_2,\ldots, e_M$ create the first $M$ eigenvectors of such a
compartmental system.

To lump the system of differential equations induced by this model,
we can choose some of these eigenvectors to generate several lumping
matrices.

In the first case, if we do not use the vector [1,1,\ldots,1,1] to
generate $Q$, only $\hat{M}$ of the first $M$ elements of the
standard basis for $\mathbb{R}^{M+1}$, that appear above, we will
receive a new compartmental system, with $\hat{M}$ compartments,
where the new species are taken from the old external ones only. In
this case lumping actually discards some peripheral compartments and
permutes the remaining ones.

If we take an $\hat{M}\times\hat{M}$ nonsingular matrix, $P,$ i.e. a
basis transformation matrix, then $PQ$ will be another lumping
matrix. It will consist of some of $P$'s columns, and values being 0
elsewhere. Accordingly the new compartments will be the linear
combinations of certain old peripheral compartments. An obvious
interpretation is that they are measured together.

Assume e.g. we have chosen $Q$ in the following way: it consists
of $e_i$, $e_j$ and $e_k$ of the natural basis
$\mathbb{R}^{M+1}$, $i,j,k\in\{1,2,\ldots, M\}$. Let $P\in \mathbb{R}^{3\times 3}$ be an invertible matrix.
Then
$$\hat{x}=PQx=x_ip_{.1}+x_jp_{.2}+x_kp_{.3},$$ where
$p_{.1},p_{.2},p_{.3}$ are the linearly independent columns of $P$, and the coordinates of the new composition vector
$\hat{x}$ are  linear combinations of the external species $x_i, x_j$
and $x_k$ with (in general) different coefficients.

Now, suppose, the eigenvector $[1,1,\ldots,1]$ is contained in the
rows of matrix $Q.$ In this case the system of
equation $\hat{x}_i=\sum_{j=1}^{M+1}q_{ij}x_j$ defines new
compartments, composed by some of the existing peripheral ones, plus
the sum of the original one.

As an example let us consider the following irreversible mamillary
system with inward flows:
\begin{eqnarray*}
&&\mathcal{X}_3\\
&&\downarrow \hbox{\scriptsize{$k_3$}}\\
\mathcal{X}_1\stackrel{k_1}{\rightarrow}&&\mathcal{X}_4\stackrel{k_2}{\leftarrow}\mathcal{X}_2
\end{eqnarray*}

The induced kinetic differential equation is $\dot{x}=Ax$, where
\[A=
 \left[
  \begin{array}{cccc}
  -k_1   & 0       &  0      & 0\\
  0      & -k_2    &  0      & 0\\
  0      & 0       &  -k_3   & 0\\
  k_1    & k_2     &  k_3    & 0
  \end{array}
 \right].
\]
Using the fact that the eigenvectors of $A^T$ are $[1,0,0,0],
[0,1,0,0], [0,0,1,0],$ and $[1,1,1,1],$ we can set e.g.
\[Q=
 \left[
  \begin{array}{cccc}
  0      & 0            &  1      & 0\\
  1      & 0            &  0      & 0\\
  1      & 1            &  1      & 1

  \end{array}
 \right].
\]
The new variables become
\[\left[
  \begin{array}{c}
  \hat{x}_1\\
  \hat{x}_2\\
  \hat{x}_3\\
  \end{array}
 \right]=\left[
  \begin{array}{cccc}
  0      & 0            &  1      & 0\\
  1      & 0            &  0      & 0\\
  1      & 1            &  1      & 1
  \end{array}
 \right]\left[
  \begin{array}{c}
  x_1\\
  x_2\\
  x_3\\
  x_4
  \end{array}
 \right].
 \] and the lumped system has the variables
\begin{equation}
\begin{split}
\hat{x}_1&=x_3\\
\hat{x}_2&=x_1\\
\hat{x}_3&=x_1+x_2+x_3+x_4.
\end{split}
\end{equation}

The resultant process obeys a differential equation
\[\left[
  \begin{array}{c}
  \dot{\hat{x}}_1\\
  \dot{\hat{x}}_2\\
  \dot{\hat{x}}_3\\
  \end{array}
 \right]=\left[
  \begin{array}{ccc}
  -k_3   & 0      &  0\\
  0      & -k_1   &  0\\
  0      & 0      & 0
  \end{array}
 \right]\left[
  \begin{array}{c}
  \hat{x}_1\\
  \hat{x}_2\\
  \hat{x}_3\\
  \end{array}
 \right]
 \]
 which is the induced kinetic differential equation e.g. of the reaction
\begin{eqnarray*} \label{example1}
\hat{\mathcal{X}}_1\stackrel{k_3}{\rightarrow}\mathcal{O}
\stackrel{k_1}{\leftarrow}\hat{\mathcal{X}}_2,
\end{eqnarray*} that is no more a mamillary system. Or, again, we can take
$\hat{\mathcal{X}}_1\rightarrow\hat{\mathcal{X}}_3
\leftarrow\hat{\mathcal{X}}_2,$ and say we are not interested in
the change of concentration of $\hat{\mathcal{X}}_3;$ we consider it as an external species.

\textbf{Outward flows}

\begin{figure}[h!]
\begin{center}
\includegraphics[totalheight=2.2in]{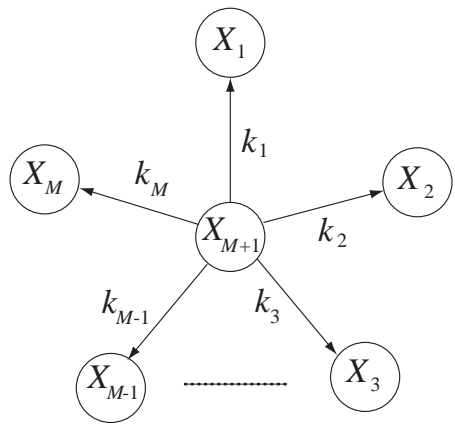}
\end{center}
\caption{Outward flows}\label{OutFlows}
\end{figure}

The induced kinetic differential equation is $\dot{x}=Ax,$ where
\[ A=
 \left[
  \begin{array}{ccccc}
  0      & 0            & \ldots & 0      & k_1\\
  0      & 0            & \ldots & 0      & k_2\\
  \vdots & \vdots       & \vdots & \vdots &\vdots\\
  0      & 0            & \ldots & 0      & k_M\\
  0      & 0            & \ldots & 0      & -K
  \end{array}
 \right],
\] with $K=(k_1+k_2+\cdots+k_M).$
The transpose of it, $A^T$, has a single eigenvalue
$-K$ with the eigenvector
$[0,0,\ldots,0,1]$, and an  eigenvalue 0 with multiplicity $M$,
with the corresponding independent eigenvectors

\[
  \begin{array}{ccccc}
  \left[ 1      \right.& 0                    & \ldots & -\frac{k_1}{k_M}                 &\left. 0\right]     \\
  \left[ 0      \right.& 1                    & \ldots & -\frac{k_2}{k_M}                 &\left. 0\right]     \\
  \vdots               & \vdots               & \vdots & \vdots                           & \vdots             \\
  \left[ 0      \right.& 0                    & \ldots & -\frac{k_{M-1}}{k_M}             &\left. 0\right]     \\
  \left[ 0      \right.& 0                    & \ldots & \frac{x}{k_M}  &\left. 1\right]
  \end{array}
 \]

\hide{
\begin{displaymath}
[1,0,\ldots,-\frac{k_1}{k_M},0],[0,1,
\ldots,-\frac{k_2}{k_M},0],\ldots, [0,0,
\ldots,-\frac{k_{M-1}}{k_M},0],
[0,0,\ldots,\frac{k_1+\cdots+k_M}{k_M},1].
\end{displaymath}
}

If we build up a lumping matrix, $Q$, we get reasonable result only
with eigenvectors belonging to the multiple eigenvalue $0$, since we
get $\hat{A}=0$ in all other cases, and it is not worth taking such
a $Q$.

If the eigenvector $[0,0,\ldots,0,1]$ appears in the lumping matrix,
we obtain a lumped system, whose coefficient matrix consists of the
$0$ elements, except a single element on the principal diagonal,
which has the value $-K.$ This can be represented by the extremely
simple reaction
$\hat{\mathcal{X}}\stackrel{K}{\rightarrow}\mathcal{O}.$

\hide{ As an example consider the following lumping matrix
\[
 \left[
  \begin{array}{ccccc}
  0  & 0  & 1  & -\frac{k_3}{k_4} & 0\\
  0  & 0  & 0  & 0                & 1\\
  0  & 0  & 0  & \frac{k_1+k_2+k_3+k_4}{k_4}      & 1
  \end{array}
 \right].
\]
}

\textbf{Irreversible mamillary systems with inward and outward flows}

Instead of giving a general treatment we shall take an example again, as in subsection \ref{nonunich}. Let us consider the mamillary system below.

\begin{figure}[h!]
\begin{center}
\includegraphics[totalheight=2.2in]{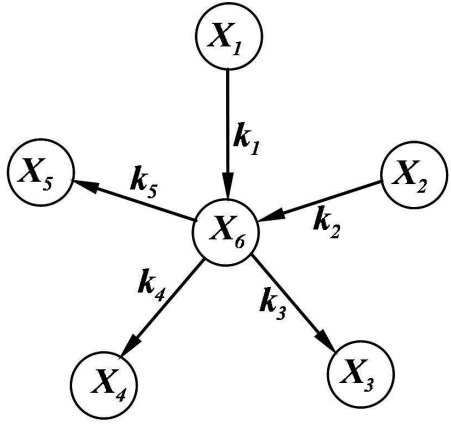}
\end{center}
\caption{An irreversible mamillary system with inward and outward flows}\label{InOutFlows}
\end{figure}

The coefficient matrix of the induced kinetic differential equation is
\[
 \left[
  \begin{array}{cccccc}
  -k_1    & 0      & 0    & 0      & 0   & 0  \\
     0    & -k_2   & 0    & 0      & 0   & 0  \\
     0    & 0      & 0    & 0      & 0   & k_3\\
     0    & 0      & 0    & 0      & 0   & k_4\\
     0    & 0      & 0    & 0      & 0   & k_5\\
     0    & 0      & 0    & 0      & 0   & -(k_3+k_4+k_5)\\
  \end{array}
 \right],
\] with the simple eigenvalues $-k_1,-k_2, -K:=-(k_3+k_4+k_5),$
and with the triple eigenvalue 0. The corresponding eigenvectors of
its transpose are $[1,0,0,0,0,0],$ $[0,1,0,0,0,0],$ $[0,0,0,0,0,1],$
$\left[0,0,1,0,0,\cfrac{k_3}{K}\right],$
$\left[0,0,1,0,-\cfrac{k_3}{k_5},0\right],$
$\left[0,0,1,-\cfrac{k_3}{k_4},0,0\right].$

Taking the lumping matrix $$Q=\left[ \begin {array}{cccccc}
0&0&1&0&0&{\cfrac {k_{{3}}}{K}}\\\noalign{\medskip}0&0&1&-{\cfrac
{k_{{3}}}{k_{{4}}}}&0&0
\\\noalign{\medskip}1&0&0&0&0&0\end {array} \right],$$ we get the
lumped system
\begin{equation}
\begin{split}
\hat{x}_1&=x_3+\cfrac{k_3}{K}x_6\\
\hat{x}_2&=x_3-\cfrac{k_3}{k_4}x_4\\
\hat{x}_3&=x_1.
\end{split}
\end{equation}

In this case the lumped system's differential equation will be very
simple:
\[\left[
  \begin{array}{c}
  \dot{\hat{x}}_1\\
  \dot{\hat{x}}_2\\
  \dot{\hat{x}}_3\\
  \end{array}
 \right]=\left[
  \begin{array}{ccc}
  0      & 0      &  0\\
  0      & 0      &  0\\
  0      & 0      &  -k_1
  \end{array}
 \right]\left[
  \begin{array}{c}
  \hat{x}_1\\
  \hat{x}_2\\
  \hat{x}_3\\
  \end{array}
 \right].
 \]
 We can associate it to the reaction
$\hat{\mathcal{X}}_3\stackrel{k_1}{\rightarrow}\mathcal{O}.$

\subsubsection{Simplicial compartmental systems}\label{classmam}
Suppose we have the following formal reaction steps as follows
\vspace{-0.5cm}
\begin{eqnarray*} \label{example1}
\mathcal{O}\stackrel{d}{\leftarrow}\mathcal{X}_i
\begin{array}{c}
\hbox{\scriptsize{$c_{j-i}$}}\\[-0.9em]
\rightleftharpoons\\[-1.4em]
\hbox{\scriptsize{$c_{M-j+i}$}}\end{array}
\mathcal{X}_j\stackrel{d}{\rightarrow}\mathcal{O},\quad \mbox{ for }
i<j;\,i,j\in\{1,2,\ldots,M\}.
\end{eqnarray*}
The fact that the reaction rate coefficients are the same for many
reaction-antireaction pairs may come from the application when the
compartments are physically separated (by a membrane e.g.) parts of
the space. In general, such kinds of assumption are made in cases
when diffusion is modeled by mass transport between homogeneous
boxes; such models often arise \cite{sh}, \cite{sherr}.

The transpose of the coefficient matrix of the induced kinetic differential equation is
\[ A^{\top}=
 \left[
  \begin{array}{ccccc}
  c_0      & c_1      & c_2  & \ldots  & c_{M-1}\\
  c_{M-1}  & c_0      & c_1  & \ldots  & c_{M-2}\\
  c_{M-2}  & c_{M-1}  & c_0  & \ldots  & c_{M-3}\\
  \vdots   &          &      &         & \vdots \\
  c_1      & c_2      & c_3  & \ldots  & c_0
  \end{array}
 \right],
\] with $c_0:=-(c_1+c_2+\cdots +c_{M-1}+d).$ Such a matrix $A^{\top}$
(for which every row is a cyclic permutation of the top row) is
called a cyclic, or circulant matrix. Its eigenvalues can be
calculated easily \cite{gray06}. (Certainly, $A$ is a cyclic matrix,
as well.)
\begin{eqnarray*}
\begin{split}
    \lambda_1&=\sum_{m=0}^{M-1}c_m\\
    \lambda_2&=\sum_{m=0}^{M-1}c_m\varepsilon_1^m\\
    \vdots\\
    \lambda_M&=\sum_{m=0}^{M-1}c_m\varepsilon_{M-1}^m\\
\end{split}
\end{eqnarray*} where $\varepsilon_k:=e^{\frac{2k\pi i}{M}},\,(k=0,1,\ldots,M-1)$
are the roots of unity. \\The corresponding eigenvectors are
\[
  \begin{array}{ccccc}
  \left[ 1      \right.& 1                    & 1                 & \ldots        &\left. 1\right]     \\
  \left[ 1      \right.& \varepsilon_1           & \varepsilon_1^2      & \ldots        &\left. \varepsilon_1^{M-1}\right]\\
  \vdots               &                      &                   &               & \vdots             \\
  \left[ 1      \right.& \varepsilon_{M-1}       & \varepsilon_{M-1}^2  & \ldots        &\left. \varepsilon_{M-1}^{M-1}\right]
  \end{array}
 \]

Here we meet a new problem which we will not discuss here further:
obviously, in the applications one needs real lumped systems.

\begin{figure}[h!]
\begin{center}
\includegraphics[totalheight=1.9in]{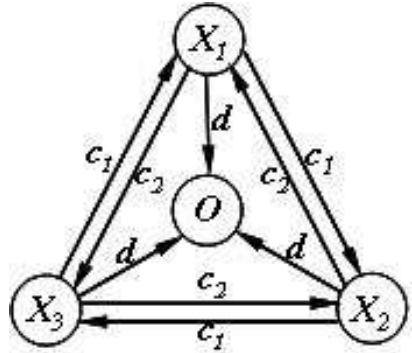}
\end{center}
\caption{A simplicial compartmental system}\label{ciklikusSik}
\end{figure}

To be more concrete, let us consider the special case (studying the
problem of complex numbers in the special case) of
Fig.~\ref{ciklikusSik}.
(This system is a special reversible circular system with outflow.)
Now
$\varepsilon_0=1,\varepsilon_1=\cfrac{-1+i\sqrt{3}}{2},\varepsilon_2=\cfrac{-1-i\sqrt{3}}{2},$
thus the eigenvalues of
\[
 \left[
  \begin{array}{ccc}
  c_0      & c_1      & c_2 \\
  c_2      & c_0      & c_1 \\
  c_1      & c_2      & c_0
  \end{array}
 \right]
\] are $\lambda_1=c_0+c_1+c_2,$ $\lambda_2=c_0+c_1\varepsilon_1+c_2
\varepsilon_1^2=c_0+c_1\cfrac{-1+i\sqrt{3}}{2}+c_2\cfrac{-1-i\sqrt{3}}{2},$
 and $\lambda_3=c_0+c_1\varepsilon_2+c_2\varepsilon_2^2=c_0+c_1
 \cfrac{-1-i\sqrt{3}}{2}+c_2\cfrac{-1+i\sqrt{3}}{2}.$

 For $c_1=c_2$ the corresponding eigenvectors are $[1,1,1], [-1,0,1], [-1,1,0].$
 In this case we can construct a few lumping matrices which lead to a new, simpler
system.

Furthermore, if $c_1\neq c_2$ we obtain the following eigenvectors:
$[1,1,1],$\\
$\left[-\cfrac{|b-c|-i\sqrt{3}(b+c)}{(2b+c)\textup{sign}(b-c)-i\sqrt{3}c},
-\cfrac{(b+2c)\textup{sign}(b-c)+i\sqrt{3}b}{(2b+c)\textup{sign}(b-c)-i\sqrt{3}c},1\right],$\\
$\left[-\cfrac{|b-c|+i\sqrt{3}(b+c)}{(2b+c)\textup{sign}(b-c)+i\sqrt{3}c},
-\cfrac{(b+2c)\textup{sign}(b-c)-i\sqrt{3}b}{(2b+c)\textup{sign}(b-c)+i\sqrt{3}c},1\right].$

\textbf{The effect of lumping on qualitative properties}

One of the major questions connected with lumping is: how are the
qualitative properties of the lumped and of the original system
connected? We investigated this problem in a more general setting in
\cite{tlrt}; here we add a new statement: suppose we lump a system
of $M$ compartments with a coefficient matrix having real
eigenvalues into a compartmental system of $\hat{M}$ compartments.
Then, none of the concentration versus time curves can have more
than $\hat{M}-2$ local extrema \cite{Pota}.

\section{A few numerical examples}\label{numer}

\subsection{A reversible chain}\label{numrevCh}
Consider a reversible chain formed by five chemical species, let
them be $\mathcal{X}_1, \mathcal{X}_2, \mathcal{X}_3, \mathcal{X}_4$
and $\mathcal{X}_5.$ Let the forward and reverse reaction rates be
$$k_1=1, k_2=4, k_3=2, k_4=1 \mbox{ and } k_{-1}=2, k_{-2}=4, k_{-3}=1,
k_{-4}=2,$$ respectively, as it can be seen in the following
chemical mechanism:

\begin{eqnarray*}
\mathcal{X}_1
\begin{array}{c}
\hbox{\scriptsize{1}}\\[-1.1em]
\rightleftharpoons\\[-1.1em]
\hbox{\scriptsize{1}}\end{array}\mathcal{X}_2
\begin{array}{c}
\hbox{\scriptsize{4}}\\[-1.1em]
\rightleftharpoons\\[-1.1em]
\hbox{\scriptsize{4}}\end{array}\mathcal{X}_3
\begin{array}{c}
\hbox{\scriptsize{2}}\\[-1.1em]
\rightleftharpoons\\[-1.1em]
\hbox{\scriptsize{5}}\end{array}\mathcal{X}_4
\begin{array}{c}
\hbox{\scriptsize{1}}\\[-1.1em]
\rightleftharpoons\\[-1.1em]
\hbox{\scriptsize{2}}\end{array}\mathcal{X}_5
\end{eqnarray*}

We can associate it with the following induced kinetic differential
equation system:{\small
$$\left\{
  \begin{array}{ll}
    \dot{x}_1=-x_1+x_2,\\
    \dot{x}_2=x_1-5x_2+4x_3,\\
    \dot{x}_3=4x_2-6x_3+5x_4,\\
    \dot{x}_4=2x_3-6x_4+2x_5,\\
    \dot{x}_5=x_4-2x_5.\\
  \end{array}
\right.$$}

The eigenvectors of the transpose of its coefficient matrix are
collected in the rows of the matrix below: {\small
$$\left[ \begin {array}{ccccc}  0.2&- 0.2&-
 0.2& 0& 1\\\noalign{\medskip}- 0.689897&
 0.069693& 0.240408& 0.449489& 1\\\noalign{\medskip}
 0.289897&- 2.869693& 4.159591&- 4.449489& 1
\\\noalign{\medskip}- 0.2& 1&- 0.2&- 2& 1
\\\noalign{\medskip} 1& 1& 1& 1& 1\end {array} \right].
$$}
If we take as a lumping matrix
$$Q=\left[ \begin {array}{ccccc}  0.289897&- 2.869693& 4.159591
&- 4.449489& 1\\\noalign{\medskip} 0.2&- 0.2&-
 0.2& 0& 1\end {array} \right],$$ after some calculations we receive
 $$\hat{A}=\left[ \begin {array}{cc} - 10.898979&
 0\\\noalign{\medskip}{ 0 }&- 2\end {array} \right]$$
and we obtain the lumped model
\begin{eqnarray*}
\hat{\mathcal{X}}_1 \stackrel{10.89}{\longrightarrow}
\hat{\mathcal{O}} \stackrel{2}{\longleftarrow} \hat{\mathcal{X}}_1.
\end{eqnarray*}

\subsection{A reversible mamillary system}\label{numrevMam}

In the following, consider a reversible compartmental system with
five compartments. Let $\mathcal{X}_5$ be the mother compartment,
and $\mathcal{X}_1, \mathcal{X}_2, \mathcal{X}_3, \mathcal{X}_4,$
the peripheral ones. Suppose that all of the reaction rates
corresponding to the reactions from the mother compartment to the
peripheral ones have the same value, $K.$ Whereas, the reverse
reactions also have identical reaction rates, $k.$

To this chemical mechanism we can set up the system $$\left\{
  \begin{array}{ll}
    \dot{x}_1= -kx_1+Kx_5\\
    \dot{x}_2= -kx_2+Kx_5\\
    \dot{x}_3= -kx_3+Kx_5\\
    \dot{x}_4= -kx_4+Kx_5\\
    \dot{x}_5= k(x_1+x_2+x_3+x_4)-4Kx_5,\\
  \end{array}
\right.$$ which describes the time evolution of the concentrations
of the species taking part in the reaction. Consequently, the
coefficient matrix will be
$$\left[ \begin {array}{ccccc} -k&0&0&0&K\\\noalign{\medskip}0&-k&0&0&K
\\\noalign{\medskip}0&0&-k&0&K\\\noalign{\medskip}0&0&0&-k&K
\\\noalign{\medskip}k&k&k&k&-4\,K\end {array} \right].
$$
Its transpose has a triple eigenvalue, $-k$, and two single
eigenvalues, 0 and and $-k-4K.$ The corresponding eigenvectors are
as follows: $[-1,0,0,1,0],$ $[-1,0,1,0,0],$ $[-1,1,0,0,0],$
$[1,1,1,1,1]$ and
$\left[-\cfrac{k}{4K},-\cfrac{k}{4K},-\cfrac{k}{4K},-\cfrac{k}{4K},1\right].$

Now, we can take several lumping matrices. For example if
$$Q=\left[ \begin {array}{ccccc} -1&0&0&1&0\\\noalign{\medskip}-
\cfrac{k}{4K}&-\cfrac {k}{4K}&-\cfrac {k}{4K}&-\cfrac
{k}{4K}&1\\\noalign{\medskip}-1&1&0&0&0\end {array} \right],
$$ we obtain the lumped system
$$\left\{
  \begin{array}{ll}
    \dot{\hat{x}}_1=-k\hat{x}_1 \\
    \dot{\hat{x}}_2=-(4K+k)\hat{x}_2\\
     \dot{\hat{x}}_3=-k\hat{x}_3,
  \end{array}
\right.$$ and the corresponding model is
\begin{eqnarray*}
\hat{\mathcal{X}}_1\stackrel{k}{\longrightarrow}&&\mathcal{O}\stackrel{k}{\longleftarrow}\hat{\mathcal{X}}_3\\
&&\uparrow \hbox{\scriptsize{$4K+k$}}\\
&&\hat{\mathcal{X}}_2
\end{eqnarray*}


\subsection{Cycles}

\subsubsection{Irreversible cycles}

\begin{figure}[h!]
\begin{center}
\includegraphics[totalheight=2.2in]{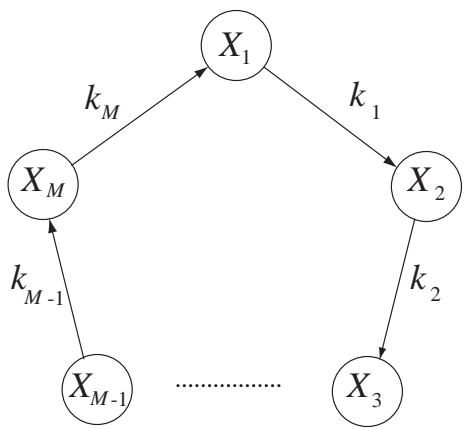}
\end{center}
\caption{Irreversible circular system}\label{IrrCirc}
\end{figure}

Consider an irreversible circular system with three compartments,
$\mathcal{X}_1,\mathcal{X}_2,\mathcal{X}_3$ and the corresponding
reaction rates $k_1,k_2,k_3.$ Then the coefficient matrix of the
induced kinetic differential equation is
$$A=\left[ \begin {array}{ccc} -k_1&0&k_3\\
\noalign{\medskip}k_1&-k_2&0\\
\noalign{\medskip}0&k_2&-k_3
\end {array} \right].$$

Since the eigenvectors of $A^{\top}$ are $[1,1,1],$\\
$\left[-\cfrac{k_1+k_2-k_3+\sqrt{k_1^2+(k_2-k_3)^2-2k_1(k_2+k_3)}}{2k_3},\right.$\\
$\left.\cfrac{k_2\big(-k_1+k_2-k_3+\sqrt{k_1^2+(k_2-k_3)^2-2k_1(k_2+k_3)}\big)}{2k_1k_3},1\right]$
and\\
$\left[\cfrac{-k_1-k_2+k_3+\sqrt{k_1^2+(k_2-k_3)^2-2k_1(k_2+k_3)}}{2k_3},\right.$\\
$\left.-\cfrac{k_2\big(k_1-k_2+k_3+\sqrt{k_1^2+(k_2-k_3)^2-2k_1(k_2+k_3)}\big)}{2k_1k_3},1\right]$
, respectively, building up $Q$ from the first two eigenvectors, in
the special case $k_1=1, k_2=2$ and $k_3=3$ we obtain the lumping
matrix
$Q=\left[ \begin {array}{ccc} 1&1&1\\
                              -\cfrac{i\sqrt{2}}{3}&\cfrac{-2+2i\sqrt{2}}{3}&1
\end {array} \right].$
After some calculations we get
$\hat{A}=\left[ \begin {array}{ccc} 0&0\\
                                    0&-3-\sqrt{2}i
\end {array} \right].$
To receive a real valued matrix $\hat{A},$ we should take, for
example $k_1=1, k_2=1/2$ and $k_3=5/128.$

\begin{figure}[h!]
\begin{center}
\includegraphics[totalheight=3.4in]{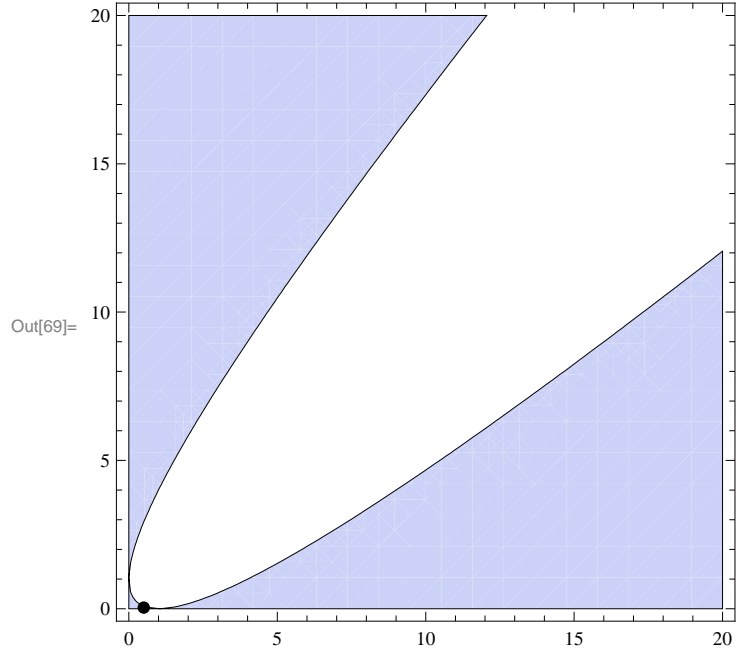}
\caption{The preferred values of $k_2$ and $k_3$ are those
\emph{outside} the curve} \label{RegionPlot}
\end{center}
\end{figure}

We can also illustrate the region of those values $k_2\in[0,20]$ and
$k_3\in[0,20],$ for which $k_1=1$ results in a lumped system with
kinetic structure, that is, a real valued matrix, $\hat{A}$ (see
fig. \ref{RegionPlot}).

\subsubsection{Reversible cycles}

Consider the reversible cycle with five compartments $\mathcal{X}_1,
\mathcal{X}_2, \mathcal{X}_3, \mathcal{X}_4, $ and $\mathcal{X}_5.$
Suppose the reaction rates are all equal to a positive real number
$k.$

Then, we can assign to this mechanism a linear differential equation
to describe the time evolution of the species's concentrations, with
coefficient matrix
$$A=\left[ \begin {array}{ccccc} -2\,k&k&0&0&k\\\noalign{\medskip}k&-2\,k
&k&0&0\\\noalign{\medskip}0&k&-2\,k&k&0\\\noalign{\medskip}0&0&k&-2\,k
&k\\\noalign{\medskip}k&0&0&k&-2\,k\end {array} \right].$$

This is a special circular matrix. $A^{\top}=A$ has two double
eigenvalues, $\cfrac{-5+\sqrt{5}}{2}$ and $\cfrac{-5-\sqrt{5}}{2},$
and a single one, 0. With the corresponding eigenvectors,
$\left[\cfrac{\sqrt{5}-1}{2},\cfrac{-\sqrt{5}+1}{2},-1,0,1\right],$
$\left[-1,\cfrac{-\sqrt{5}+1}{2},\cfrac{\sqrt{5}-1}{2},1,0\right],$\\
$\left[\cfrac{-\sqrt{5}-1}{2},\cfrac{\sqrt{5}+1}{2},-1,0,1\right],$
$\left[1,\cfrac{\sqrt{5}+1}{2},\cfrac{-\sqrt{5}-1}{2},1,0\right],$
and $[1,1,1,1,1]$ we can determine several invariant subspaces in
order to find lumping matrices. Choose, for example,
$$Q=\left[ \begin {array}{ccccc} 1&1&1&1&1\\\noalign{\medskip}-1&
\frac{\sqrt{5}+1}{2}&\frac{-\sqrt{5}-1}{2}&1&0\\\noalign{\medskip}
\frac{-\sqrt{5}-1}{2}&\frac{\sqrt{5}+1}{2}&-1&0&1\end {array}
\right].$$ In this case the lumped system will be
\[\left[
  \begin{array}{c}
  \dot{\hat{x}}_1\\
  \dot{\hat{x}}_2\\
  \dot{\hat{x}}_3\\
  \end{array}
 \right]=\left[
  \begin{array}{ccc}
  0      & 0      &  0\\
  0      & -\frac{5+\sqrt{5}}{2}k      &  0\\
  0      & 0      & -\frac{5+\sqrt{5}}{2}k
  \end{array}
 \right]\left[
  \begin{array}{c}
  \hat{x}_1\\
  \hat{x}_2\\
  \hat{x}_3\\
  \end{array}
 \right].
 \]
 We can associate it to the model
$\hat{\mathcal{X}}_2\stackrel{\frac{5+\sqrt{5}}{2}k}
{\longrightarrow}\mathcal{O}\stackrel{\frac{5+\sqrt{5}}{2}k}
{\longleftarrow}\hat{\mathcal{X}}_3.$

\section{Discussion, plans}
The most important classes of compartmental systems have been
reviewed from the point of view of symbolic lumpability. Practically
interesting lumped systems mainly arise from numerical calculations,
which can be carried out in all cases without difficulties. We used
the sentence "which is the induced kinetic differential equation of
the reaction" recurrently. However, given a kinetic differential
equation the inducing reaction is by far not unique\\ \cite[pages
67--69]{ET}.

\section{Appendix}

\textit{Suppose we are given two natural numbers, $n$ and $\hat{n}$,
$\hat{n}\leq n$, and an $\hat{n}\times n$ matrix $Q$ of full rank
with real elements. The question arises: what are the necessary and
sufficient conditions for the existence of a nonsingular
$\hat{n}\times\hat{n}$ matrix $P$ such that all elements of $PQ$ are
nonnegative?}

This question is hard enough to answer in a general case. Here is a
result when $\hat{n}\leq2.$ One can see, if $\hat{n}=1$, the
elements of $Q$ must have identical sign, for the existence of such
$P.$

Now, assume $\hat{n}=2,$ and take
\begin{equation}\label{Q}
Q=
 \left[
  \begin{array}{cccc}
  a_{11}      & a_{12}    & \ldots  & a_{1n}\\
  a_{21}      & a_{22}    & \ldots  & a_{2n}
  \end{array}
 \right].
\end{equation}
Furthermore, take
\begin{equation}\label{P}
P=
 \left[
  \begin{array}{cc}
  p_{11}      & p_{12}\\
  p_{21}      & p_{22}
  \end{array}
 \right].
\end{equation}
Then we obtain
\[PQ=
 \left[
  \begin{array}{cccc}
  a_{11}p_{11}+a_{21}p_{12}      & a_{12}p_{11}+a_{22}p_{12}    & \ldots  & a_{1n}p_{11}+a_{2n}p_{12}\\
  a_{11}p_{21}+a_{21}p_{22}      & a_{12}p_{21}+a_{22}p_{22}    & \ldots  & a_{1n}p_{21}+a_{2n}p_{22}
  \end{array}
 \right].
\]
The requirement is that all the elements of the matrix above should
be nonnegative real numbers. This assumption, i.e. the inequalities
$a_{1j}p_{i1}+a_{2j}p_{i2}\geq0,\,\,j=1,\ldots,n$ determine
half-planes in the plane $(p_{i1},p_{i2}),$ passing through the
origin, $i=1,2.$ Thus, the problem is to find the cases, when the
intersections of the corresponding planes (which is in accordance
with the first, resp. the second row in $PQ$) are not empty.

Examine the columns of the matrix $Q.$ We distinguish 9 cases. In
what follows, the symbols $+$, and $-$ indicate the presence of a
positive or a negative number in the matrix $Q$.

\begin{tabular}{ccccccccc}
1.&2.&3.&4.&5.&6.&7.&8.&9.\\
\(\left[\begin{array}{c}0\\0\end{array}\right]\)&
\(\left[\begin{array}{c}0\\-\end{array}\right]\)&
\(\left[\begin{array}{c}+\\-\end{array}\right]\)&
\(\left[\begin{array}{c}+\\0\end{array}\right]\)&
\(\left[\begin{array}{c}+\\+\end{array}\right]\)&
\(\left[\begin{array}{c}0\\+\end{array}\right]\)&
\(\left[\begin{array}{c}-\\+\end{array}\right]\)&
\(\left[\begin{array}{c}-\\0\end{array}\right]\)&
\(\left[\begin{array}{c}-\\-\end{array}\right]\)
\end{tabular}

\begin{figure}[h!]
\begin{center}
\includegraphics[totalheight=3.4in]{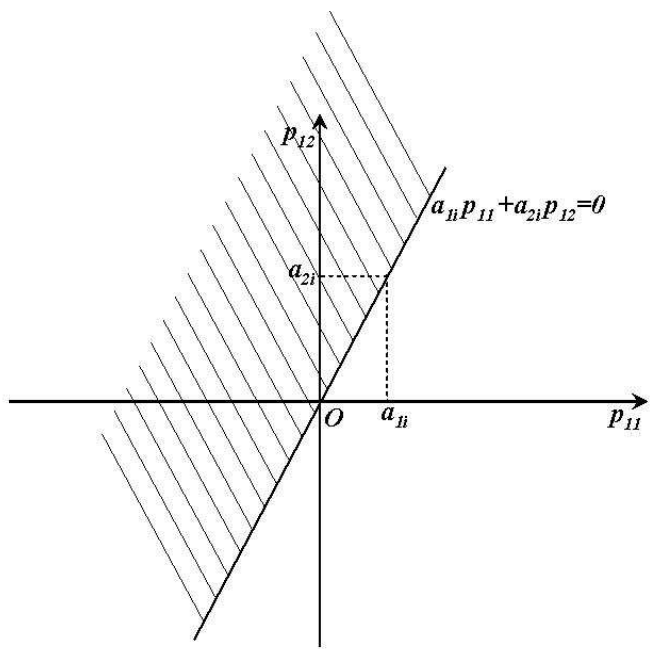}
\caption{A graphical representation for the case 3.}\label{felsik}
\end{center}
\end{figure}

For example, suppose the $i$th column of $Q$ is of type 3, i.e.
$a_{1i}$ is positive and $a_{2i}$ is negative. Then the inequality
$a_{1i}p_{11}+a_{2i}p_{12}\geq0$ corresponds to the case in the
fig.~\ref{felsik}.

The slope of the line with equation $a_{1i}p_{11}+a_{2i}p_{12}=0$
depend on the number $-a_{1i}/a_{2i}.$  Here this fraction is a
positive number. The shaded region in the figure represents the
region that is excluded from the solution.

Thus, after the geometrical consideration, we can conclude: if the
matrix $Q$ contains columns of form $\left[\begin{smallmatrix}
                        a \\
                        b
       \end{smallmatrix}
  \right]$ and
$\left[\begin{smallmatrix}
                        -a \\
                        -b
       \end{smallmatrix}
  \right]$ simultaneously, $\forall \,\, a,b\in\mathbb{R},$ with $a$ and $b$
  different from 0 at the same time, there does not exist nonsingular,
  $2\times2$ matrix $P,$ which satisfies the requirement $PQ\geq0.$
  This is the case when matrix $Q$ has a pair of columns of type 2 and 6, or 4 and 8, or 3 and 7, or 5 and 9. In the last two cases the elements could only differ in sign. Henceforth, for a shortest notation we will use 26,48,37,59 to point to pair of cases when the matrix $P$ does not exist.

Finding all the cases when the matrix $Q$ contains three columns
which precisely exclude together the existence of $P,$ we lean on
the geometrical representation again. Assume $Q$ does not contain a
pair of columns fitting the case described above. If we check the
three half-plane cases, we get the following result: 247, 257, 258,
358, 368, 369, 469, 479 and  569 are the cases that exclude each
other, i.e. the existence of $P,$ by all means. Furthermore, there
are other instances for the nonexistence of such a $P$, 259, 347,
357, 359, 367, 378, 379 and 459, but in these cases we still have to
verify another condition regarding the slopes. This fact will be
illustrated later in an example.

In the case of the intersection of four half-planes, assumed that we
did not find in $M$ columns corresponding to either cases given
earlier, we get only one case for empty intersection, specifically
for 2358, i.e. for a matrix that contains the columns
$\left[\begin{smallmatrix} 0 \\ - \end{smallmatrix}  \right],$
$\left[\begin{smallmatrix} + \\ - \end{smallmatrix}  \right],$
$\left[\begin{smallmatrix} + \\ + \end{smallmatrix}  \right],$
$\left[\begin{smallmatrix} - \\ 0 \end{smallmatrix}  \right] $
together, disregarding the order.

The cases presented previously exhaust all the cases, when to a
$2\times n$ matrix given in \eqref{Q} we cannot find a $2\times 2$
nonsingular matrix, $P,$ so that all the elements of $PQ$ are
nonnegative.

Let us take an example to illustrate the problem. Consider the matrix
$$Q=\left[
  \begin{array}{rrrr}
   5   & 2   & 2  & -3\\
  -2   & 0   & 1  & -1
  \end{array}
 \right],$$ and $P$ as in \eqref{P}. Then
 $$PQ=\left[
  \begin{array}{rrrr}
   5p_{11}-2p_{12}   & 2p_{11}   & 2p_{11}+p_{12}  & -3p_{11}-p_{12}\\
   5p_{21}-2p_{22}   & 2p_{21}   & 2p_{21}+p_{22}  & -3p_{21}-p_{22}
  \end{array}
 \right],$$ whose elements must satisfy the system of inequalities:
   \begin{equation}\label{ineq}
    \begin{cases}
    5p_{11}-2p_{12}& \geq  0, \\
            2p_{11}& \geq  0, \\
     2p_{11}+p_{12}& \geq  0, \\
    -3p_{11}-p_{12}& \geq  0, \\
    5p_{21}-2p_{22}& \geq  0, \\
            2p_{21}& \geq  0, \\
     2p_{21}+p_{22}& \geq  0, \\
    -3p_{21}-p_{22}& \geq  0.
    \end{cases}
\end{equation}
Consider the first four of them, and give a geometrical representation as in figure~\ref{rossz}.

\hide{
\begin{figure}[h!]
\begin{center}
\includegraphics[totalheight=3.0in]{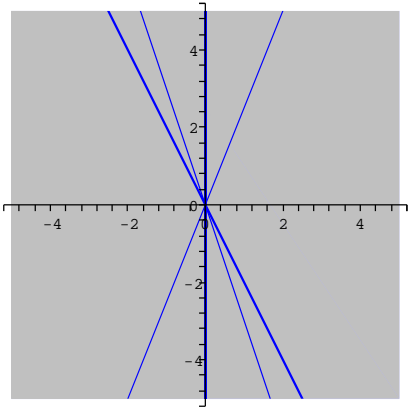}
\vspace{-1.0cm}
\caption{Example for nonexistence of $P$}\label{rossz}
\end{center}
\end{figure}
}

\begin{figure}[!h]
    \begin{minipage}[b]{0.5\linewidth} 
        \centering
        \includegraphics[width=6.7cm]{pelda01.eps}
        \vspace{-1.0cm}
        \caption{Example for nonexistence of $P$}\label{rossz}
    \end{minipage}
        \hspace{0.2cm} 
    \begin{minipage}[b]{0.5\linewidth}
        \centering
        \includegraphics[width=6.7cm]{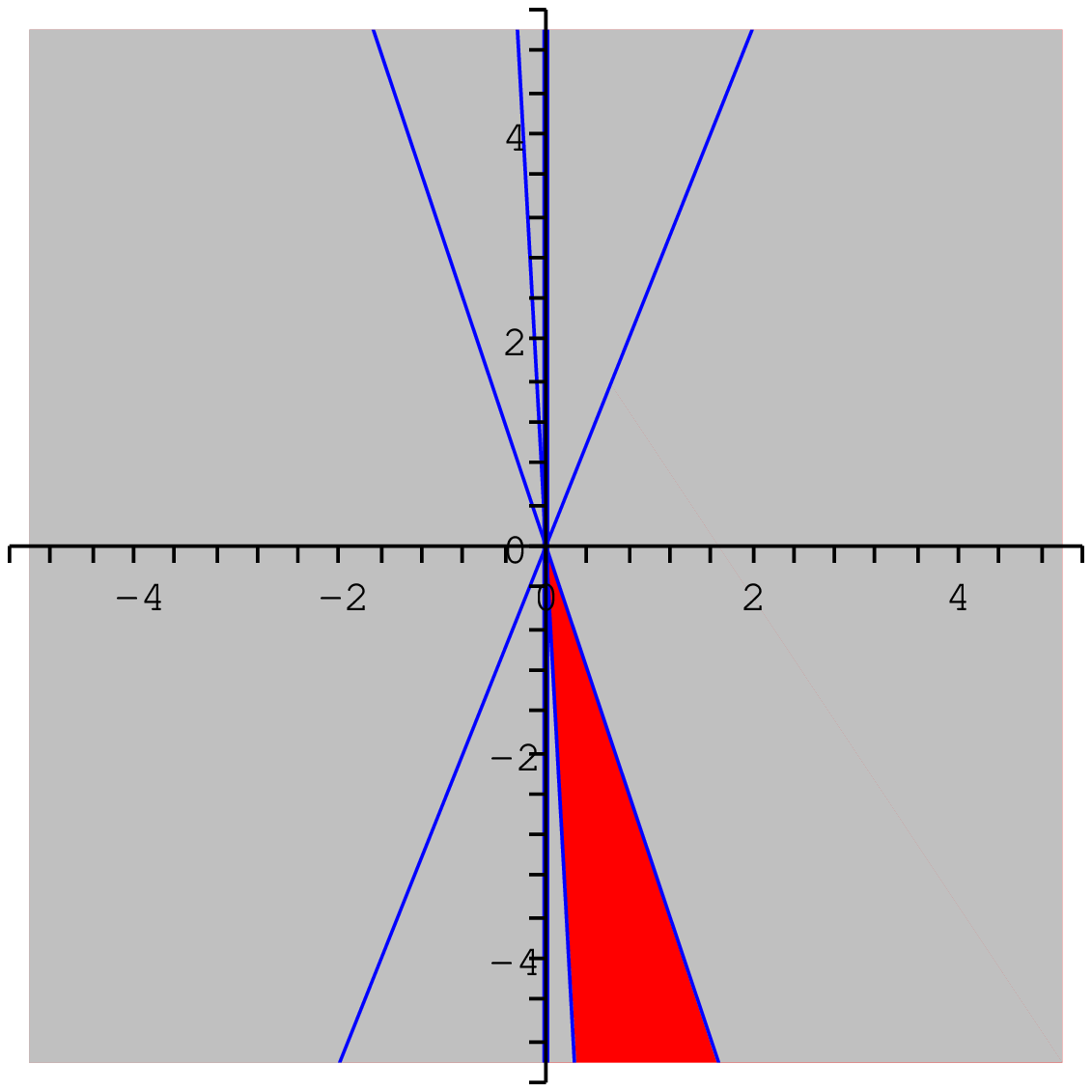}
        \vspace{-1.0cm}
        \caption{Example for existence of $P$}\label{jo}
    \end{minipage}
\end{figure}

The picture was created with Maple program, and one can see, that in
this case the inequality system \eqref{ineq} does not have any
solution. This was only to be expected, because the matrix contains
columns of type 3, 4, 5 and 9, and 359 is a critical case, since the
slope of the line corresponding to case 9 is smaller then the one
corresponding to 5, i.e. $-(-3)/(-1)<-2/1.$

Notice that if we take $a_{13}=18$ in $Q,$ the situation will
change, i.e. the existence of $P$ will be insured, because in this
case the direction of the inequality regarding the slopes will
change, as it can be seen in Fig.~\ref{jo}.

\hide{
\begin{figure}[h!]
\begin{center}
\includegraphics[totalheight=3.0in]{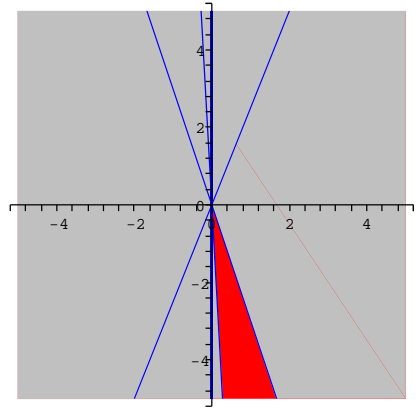}
\vspace{-1.0cm}
\caption{Example for existence of $P$}\label{jo}
\end{center}
\end{figure}
}

Choose a point from the region that indicates the solution (the dark
region in Fig. \ref{jo}), such as $(p_{11},p_{12})=(1,-5).$ To
determine $p_{21}$ and $p_{22}$ one must choose their values from
the same region, except the case when the point $(p_{21},p_{22})$
can be found on the line defined by the origin and the point
$(p_{11},p_{12})=(1,-5).$ The justification of this statement is as
follows. The equation of the line that passes through the points
$(0,0)$ and $(p_{11},p_{12})$ is $y=(p_{12}/p_{11})x.$ If
$(p_{21},p_{22})$ is a point on this line, then it must satisfy the
relation $p_{11}p_{22}-p_{12}p_{21}=0,$ which is equivalent to the
condition $\det{P}=0,$ but $P$ cannot be singular. For example we
can choose $(p_{21},p_{22})=(1/2,-4).$ In this case $\det P=-3/2.$
Now, check the nonnegativity of the elements of the matrix $PQ$:
$$PQ=
 \left[
  \begin{array}{cc}
  1           & -5   \\
  \frac{1}{2} & -4
  \end{array}
 \right]
 \left[
  \begin{array}{rrrr}
   5   & 2   & 18  & -3\\
  -2   & 0   & 1  & -1
  \end{array}
 \right]=
 \left[
  \begin{array}{rrrr}
   15         & 2   & 13  & 2\\
   \frac{21}{2}  & 1   & 5   & \frac{5}{2}
  \end{array}
 \right].
$$

\textit{Suppose we are given two natural numbers, $n$ and $\hat{n}$, $\hat{n}\leq n$, and an
$\hat{n}\times n$ matrix $Q$ of full rank with complex elements. The question arises, what are the
necessary and sufficient conditions for the
existence of a nonsingular $\hat{n}\times\hat{n}$ matrix $P$ such that all elements
of $PQ$ are real?}

Observe that for $\hat{n}=1,$ $Q$ must have a special form to find a suitable $P$ to it, i.e. we must have either
$Q=\left[im_1\,\, im_2\, \ldots\, im_n\right],$ with $m_j\in \mathbb{R},\, j=1,\ldots,n,$  or all the elements of $Q$ have to be real, but this is a trivial case.

Consider $\hat{n}>1$ and $Q$ an $n\times\hat{n}$ matrix, with $Q=Q_1+iQ_2$ where $Q_1,Q_2\in \mathcal{M}_{\hat{n}n}(\mathbb{R}).$ A sufficient condition for the existence of a $\hat{n}\times\hat{n}$ matrix $P$ such that all elements of $PQ$ are real is:
\begin{enumerate}
\item $Q_1$ and $Q_2$ are nonsingular;
\item $Q_1^TQ_2=Q_2^TQ_1.$
\end{enumerate}

In this case we can choose $P=Q_1^T-iQ_2^T.$ Realize that if the
conditions above are satisfied, then {\setlength\arraycolsep{2pt}
\begin{eqnarray*}
    PQ & = &(Q_1^T-iQ_2^T)(Q_1+iQ_2)= Q_1^TQ_1+iQ_1^TQ_2-iQ_2^TQ_1+Q_2^TQ_2\\
    & = & Q_1^2+Q_2^2+i(Q_1^TQ_2-Q_2^TQ_1)\in\mathcal{M}_{\hat{n}n}(\mathbb{R}).
\end{eqnarray*}
}

In what follows suppose that the conditions above are not satisfied
and take
$$Q=
 \left[
  \begin{array}{cccc}
 a_{11}+ib_{11} & a_{12}+ib_{12} & \ldots  & a_{1n}+ib_{1n}\\
 a_{21}+ib_{21} & a_{22}+ib_{22} & \ldots  & a_{2n}+ib_{2n}\\
 \vdots & \vdots & \vdots & \vdots \\
 a_{\hat{n}1}+ib_{\hat{n}1} & a_{\hat{n}2}+ib_{\hat{n}2} & \ldots  & a_{\hat{n}n}+ib_{\hat{n}n}
  \end{array}
 \right],
$$
$$
P=
 \left[
  \begin{array}{cccc}
  p_{11}+iq_{11} & p_{12}+iq_{12} & \ldots  & p_{1\hat{n}}+iq_{1\hat{n}}\\
  p_{21}+iq_{21} & p_{22}+iq_{22} & \ldots  & p_{2\hat{n}}+iq_{2\hat{n}}\\
  \vdots & \vdots & \vdots & \vdots \\
 p_{\hat{n}1}+iq_{\hat{n}1} & p_{\hat{n}2}+iq_{\hat{n}2} & \ldots  & p_{\hat{n}\hat{n}}+iq_{\hat{n}\hat{n}}
  \end{array}
 \right].
$$

The requirement that all the elements of $PQ$ are real is equivalent
to a linear, homogeneous system
\begin{equation}\label{nagyrendsz}
\sum_{j=1}^{\hat{n}}(p_{lj}b_{jk}+q_{lj}a_{jk})=0,\quad k=\overline{1,n},\,\, l=\overline{1,\hat{n}},
\end{equation}
where the number of unknowns are $2\hat{n}^2$ and the number of
equations are $\hat{n}n.$ Notice that this system has an interesting
property: it can be divided into independent subsystems regarding
the unknowns. In this case each of these subsystems can be solved
separately. Moreover, all of the results can be written in identical
form, since they have the same coefficient-matrix. Therefore, taking
$l=1$, consider and treat only the subsystem with $n$ equations
\begin{equation}\label{rendsz}
\sum_{j=1}^{\hat{n}}(p_{1j}b_{jk}+q_{1j}a_{jk})=0,\quad k=\overline{1,n},
\end{equation} with unknowns $p_{11}, p_{12},\ldots, p_{1,\hat{n}},$ and $q_{11}, q_{12},\ldots, q_{1,\hat{n}}.$
Depending on $n$ and $\hat{n}$ the number of unknowns can be smaller or bigger than the number of equations.
The coefficient-matrix is
\begin{equation}
\left[
  \begin{array}{cccccccc}\label{matr}
  b_{11} & b_{21} & \ldots  & b_{\hat{n}1} & a_{11} & a_{21} & \ldots  & a_{\hat{n}1}\\
  b_{12} & b_{22} & \ldots  & b_{\hat{n}2} & a_{12} & a_{22} & \ldots  & a_{\hat{n}2}\\
  \vdots & \vdots & \vdots & \vdots & \vdots & \vdots & \vdots & \vdots\\
  b_{1n} & b_{2n} & \ldots  & b_{\hat{n}n} & a_{1n} & a_{2n} & \ldots  & a_{\hat{n}n}\\
  \end{array}
 \right].
\end{equation}

Denote by $d$ the main determinant for \eqref{matr}. Then, if $d\neq0,$ system \eqref{rendsz} has identical zero solution: $\underbrace{(0,\ldots, 0)}_{2\hat{n}}.$ In this case $\det(P)=0,$ so, $P$ does not satisfy the nonsingularity requirement.

For an adequate $P$ we must have $d=0.$ If it holds, we must specify
the rank of \eqref{matr}. Let us denote it by $r.$ If $r=2\hat{n},$
i.e. the system \eqref{rendsz} is determined, and has a unique
solution. But in this case all the other subsystems mentioned above
for \eqref{nagyrendsz} has the same (constant) solution. Thus, the
matrix $P$ will consist of identical rows, so we have $\det(P)=0$
again.

For $r=2\hat{n}-1$ the system is indefinite, and the solution has
the form
$$(\alpha_{1}c_{1}, \alpha_{2}c_{2},\ldots, \alpha_{2\hat{n}}c_{2\hat{n}}),$$ where $\alpha_{i}$ are parameters and $c_i$ represent constants, $i=1,\ldots,2\hat{n}.$
Consequently, $P$ can be given in such a way that its rows will be $k$-times for the others. Thus, we discover a nonsingular $P$ once more.

Finally, if $r<2\hat{n}-1,$ then there exists an adequate matrix
$P,$ because in this case the solution of \eqref{rendsz} can be
expressed with at least two parameters, and this gives possibility
for choosing linearly independent lines to $P.$

Now, look at an example. Consider
$$Q=\left[
  \begin{array}{cccc}
   1+i   & 2+i   & 4+2i  & 2+2i\\
  -1     & 2i    & 4i    & -2
  \end{array}
 \right], \quad \mbox{and take} \quad
P=\left[
  \begin{array}{cc}
   p_{11}+iq_{11}   & p_{12}+iq_{12}\\
   p_{21}+iq_{21}   & p_{22}+iq_{22}
  \end{array}
 \right].
 $$

In this case the coefficient-matrix \eqref{matr} is
$$\left[
  \begin{array}{cccc}
   1 & 0 & 1 & -1\\
   1 & 2 & 2 &  0\\
   2 & 4 & 4 &  0\\
   2 & 0 & 2 & -2
  \end{array}
 \right],
$$ and for it $r=2<2\hat{n}-1.$ Thus, we can compute a nonsingular
$P\in\mathcal{M}_{22}(\mathbb{C}).$ Since the homogeneous, linear
system with coefficient-matrix above has solution $(-\alpha+\beta,
-\frac{1}{2}(\alpha+\beta),\alpha, \beta),$ we can take $p_{11}=0,
p_{12}=1, q_{11}=-1, q_{12}=-1$ and similarly $p_{21}=2, p_{22}=2,
q_{21}=-3, q_{22}=-1.$ After some calculations we obtain
$$PQ=\left[
  \begin{array}{cccc}
   0 & 3 & 6  & 0\\
   3 & 9 & 18 & 8
   \end{array}
 \right].$$

\end{document}